\documentclass{article}
\usepackage{comment}
\usepackage{xcolor,colortbl}
\usepackage{graphicx} 
\usepackage{color}
\usepackage[hidelinks]{hyperref}
\usepackage{amsmath}
\usepackage{ stmaryrd }
\usepackage{amssymb}
\usepackage{amsthm}
\usepackage{amsmath}
\usepackage{mathtools}
\usepackage{fullpage}
\usepackage{expl3}

\newcommand{\pending}{\textnormal{---}}

\newtheorem{theorem}{Theorem}[section]

\newcommand{\vc}[1]{\ensuremath{\vcenter{\hbox{#1}}}}

\usepackage{tikz}
\usetikzlibrary{shapes.geometric}
\usetikzlibrary{shapes.multipart}
\usetikzlibrary{arrows}
\usetikzlibrary{calc}

\tikzset{vtx/.style={inner sep=1.7pt, outer sep=0pt, circle, fill,draw}}

\tikzset{vtxB/.style={inner sep=1.7pt, outer sep=0pt, circle, fill=white,draw}}

\tikzset{flag_pic/.style={scale=1}}  
\tikzset{unlabeled_vertex/.style={inner sep=1.7pt, outer sep=0pt, circle, fill}}
\tikzset{labeled_vertex/.style={inner sep=3pt, outer sep=0pt, rectangle, fill=white, draw=black}}
\tikzset{edge_color0/.style={color=black,line width=1.2pt,opacity=0.5,dashed}}
\tikzset{edge_color1/.style={color=red,  line width=1.2pt,opacity=0}} 
\tikzset{edge_color2/.style={color=blue, line width=1.2pt,opacity=1}}
\tikzset{edge_color3/.style={color=green!80!black,line width=1.5pt,densely dotted}}
\tikzset{edge_color4/.style={color=red, line width=1pt,decorate,decoration={zigzag,segment length=1mm,amplitude=0.2mm},line join=bevel}}
\tikzset{edge_color5/.style={color=red,  line width=1.2pt,dotted}}
\tikzset{edge_color6/.style={color=blue, line width=1.2pt,dotted}}
\tikzset{edge_color7/.style={color=green!60!black, line width=1.2pt,dotted}}
\tikzset{edge_color8/.style={color=gray, line width=1.2pt}}
\tikzset{edge_color9/.style={color=gray, dotted, line width=1.2pt}}
\tikzset{edge_color10/.style={color=gray, dashed, line width=1.2pt}}
\tikzset{edge_color11/.style={color=pink, dashed, line width=1.2pt}}
\tikzset{edge_colorroot/.style={color=red, line width=1.7pt}}
\tikzset{edge_thin/.style={color=black}}
\tikzset{edge_hidden/.style={color=black,dotted,opacity=0}}
\tikzset{vertex_color0/.style={inner sep=1.7pt, outer sep=0pt, draw, circle, fill=black}}
\tikzset{vertex_color1/.style={inner sep=1.7pt, outer sep=0pt, draw, circle, fill=red}}
\tikzset{vertex_color2/.style={inner sep=1.7pt, outer sep=0pt, draw, circle, fill=blue}}
\tikzset{vertex_color3/.style={inner sep=1.7pt, outer sep=0pt, draw, circle, fill=green!80!black}}
\tikzset{vertex_color4/.style={inner sep=1.7pt, outer sep=0pt, draw, circle, fill=pink}}
\tikzset{vertex_color5/.style={inner sep=1.7pt, outer sep=0pt, draw, circle, fill=gray,label=below:{$5$}}}
\tikzset{vertex_color6/.style={inner sep=1.7pt, outer sep=0pt, draw, circle, fill=gray,label=below:{$6$}}}
\tikzset{vertex_color7/.style={inner sep=1.7pt, outer sep=0pt, draw, circle, fill=gray,label=below:{$7$}}}
\tikzset{vertex_color8/.style={inner sep=1.7pt, outer sep=0pt, draw, circle, fill=gray,label=below:{$8$}}}
\tikzset{vertex_color9/.style={inner sep=1.7pt, outer sep=0pt, draw, circle, fill=gray,label=below:{$9$}}}
\tikzset{vertex_color10/.style={inner sep=1.7pt, outer sep=0pt, draw, circle, fill=gray,label=below:{$10$}}}
\tikzset{vertex_color11/.style={inner sep=1.7pt, outer sep=0pt, draw, circle, fill=gray,label=below:{$11$}}}
\tikzset{vertex_color12/.style={inner sep=1.7pt, outer sep=0pt, draw, circle, fill=gray,label=below:{$12$}}}
\tikzset{vertex_color13/.style={inner sep=1.7pt, outer sep=0pt, draw, circle, fill=gray,label=below:{$13$}}}
\tikzset{vertex_color14/.style={inner sep=1.7pt, outer sep=0pt, draw, circle, fill=gray,label=below:{$14$}}}
\tikzset{labeled_vertex_color0/.style={inner sep=2.2pt, outer sep=0pt, draw, rectangle, fill=black}}
\tikzset{labeled_vertex_color1/.style={inner sep=2.2pt, outer sep=0pt, draw, rectangle, fill=red}}
\tikzset{labeled_vertex_color2/.style={inner sep=2.2pt, outer sep=0pt, draw, rectangle, fill=blue}}
\tikzset{labeled_vertex_color3/.style={inner sep=2.2pt, outer sep=0pt, draw, rectangle, fill=green}}
\tikzset{labeled_vertex_color4/.style={inner sep=2.2pt, outer sep=0pt, draw, rectangle, fill=pink}}
\tikzset{labeled_vertex_color5/.style={inner sep=2.2pt, outer sep=0pt, draw, rectangle, fill=gray,label=below:{$5$}}}
\tikzset{labeled_vertex_color6/.style={inner sep=2.2pt, outer sep=0pt, draw, rectangle, fill=gray,label=below:{$6$}}}
\tikzset{labeled_vertex_color7/.style={inner sep=2.2pt, outer sep=0pt, draw, rectangle, fill=gray,label=below:{$7$}}}
\tikzset{labeled_vertex_color8/.style={inner sep=2.2pt, outer sep=0pt, draw, rectangle, fill=gray,label=below:{$8$}}}
\tikzset{labeled_vertex_color9/.style={inner sep=2.2pt, outer sep=0pt, draw, rectangle, fill=gray,label=below:{$9$}}}
\tikzset{text_color0/.style={color=black}}
\tikzset{text_color1/.style={color=red}}
\tikzset{text_color2/.style={color=blue}}
\tikzset{text_color3/.style={color=green!70!black}}
\tikzset{text_color4/.style={color=orange}}
\tikzset{text_color5/.style={color=gray}}

\def\outercycle#1#2{
\pgfmathtruncatemacro{\plusone}{#1+1}
\pgfmathtruncatemacro{\zeroshift}{270 - (#2-1)*360/#1/2 }
    \draw  \foreach \x in {0,1,...,#1}{(\zeroshift+\x*360/#1:1) coordinate(x\x)};
}

\def\outercycleH#1#2{
\pgfmathtruncatemacro{\plusone}{#1+1}
\pgfmathtruncatemacro{\minusone}{#1-1} 
\draw  \foreach \x in {0,...,\minusone}{(0.5*\x,0) coordinate(x\x)};}

\def\drawhypervertex#1#2{ \pgfmathtruncatemacro{\plusone}{#1+1}  \draw[edge_color2] (x#1)++(0,-0.2-0.2*#2)+(-0.2,0) -- +(0.2,0) +(0,0) node[fill=white,outer sep=0,inner sep=0]{\tiny \plusone};}

\def\drawhyperedge#1#2{\pgfmathtruncatemacro{\plusone}{#1+1} \draw[dotted] (x0)++(0,-0.2-0.2*#1)--++(0.5*#2-0.5,0);}

\def\labelvertex#1{\pgfmathtruncatemacro{\vertexlabel}{#1+1 } \draw (x#1) node{\color{black}\tiny\vertexlabel}; }

\newcommand{\drawHEthree}[5]{
\drawhyperedge{#2}{#1}
\drawhypervertex{#3}{#2}
\drawhypervertex{#4}{#2}
\drawhypervertex{#5}{#2}
}

\tikzset{
hyperedge/.style={fill=pink,opacity=0.5,draw=black},
}

\usepackage{ifthen}

\tikzset{vertex_u/.style={unlabeled_vertex}}
\tikzset{vertex_l/.style={labeled_vertex}}

\newcommand{\Fuu}[1]{
\,\vc{\begin{tikzpicture}[scale=0.3]\outercycle{2}{1}
\draw[edge_color#1] (x0)--(x1);
\draw (x0) node[unlabeled_vertex]{};\draw (x1) node[unlabeled_vertex]{};
\end{tikzpicture}}
\,
}

\newcommand{\Flu}[1]{
\,\vc{\begin{tikzpicture}[scale=0.3]\outercycle{2}{1}
\draw[edge_color#1] (x0)--(x1);
\draw (x0) node[labeled_vertex]{};\draw (x1) node[unlabeled_vertex]{};
\labelvertex0
\end{tikzpicture}}
\,
}

\newcommand{\Fuuu}[3]{
\vc{\begin{tikzpicture}[scale=0.4]\outercycle{3}{1}
\draw[edge_color#1] (x0)--(x1);\draw[edge_color#2] (x0)--(x2);  \draw[edge_color#3] (x1)--(x2);
\draw (x0) node[unlabeled_vertex]{};\draw (x1) node[unlabeled_vertex]{};\draw (x2) node[unlabeled_vertex]{};
\end{tikzpicture}}}

\newcommand{\Fluu}[3]{
\vc{\begin{tikzpicture}[scale=0.4]\outercycle{3}{1}
\draw[edge_color#1] (x0)--(x1);\draw[edge_color#2] (x0)--(x2);  \draw[edge_color#3] (x1)--(x2);
\draw (x0) node[labeled_vertex]{};\draw (x1) node[unlabeled_vertex]{};\draw (x2) node[unlabeled_vertex]{};
\labelvertex0
\end{tikzpicture}}}

\newcommand{\FfourEdges}[6]{
\draw[edge_color#1] (x0)--(x1);\draw[edge_color#2] (x0)--(x2);\draw[edge_color#3] (x0)--(x3);  \draw[edge_color#4] (x1)--(x2);\draw[edge_color#5] (x1)--(x3);  \draw[edge_color#6] (x2)--(x3);
}
\newcommand{\Ffour}[5]{
\vc{\begin{tikzpicture}[scale=0.4]\outercycle{4}{2}
\FfourEdges#5
\draw (x0) node[vertex_#1]{};\draw (x1) node[vertex_#2]{};\draw (x2) node[vertex_#3]{};\draw (x3) node[vertex_#4]{};
\ifthenelse{\equal{#1}{l}}{\labelvertex{0}}{}%
\ifthenelse{\equal{#2}{l}}{\labelvertex{1}}{}%
\ifthenelse{\equal{#3}{l}}{\labelvertex{2}}{}%
\ifthenelse{\equal{#4}{l}}{\labelvertex{3}}{}%
\end{tikzpicture}}
}

\usepackage{listofitems}
\newcommand{\Fnv}[2]{ 
\ifnum#2<#1 \draw (x#2) node[vertex_l]{}; \labelvertex{#2}  
\else  \draw (x#2) node[vertex_u]{}; \fi 
}
\newcommand{\Fne}[3]{
\draw[edge_color#3] (x#1)--(x#2); 
}
\newcounter{Fneid} 
\newcommand{\Fn}[3]{
\ifnum#1=1
  \vc{\begin{tikzpicture}[scale=0.4]\outercycle{1}{2}
  \Fnv{#2}{0}
  \end{tikzpicture}}
\else
\setsepchar{ }
\readlist\elabel{#3}
\pgfmathtruncatemacro{\vertexloop}{#1-1}
\pgfmathtruncatemacro{\vertexloopi}{#1-2}
\pgfmathtruncatemacro{\expectededges}{#1*(#1-1)/2}
\ifnum\elabellen=\expectededges%
\def\cycleshift{2}
\ifnum#1=2\def\cycleshift{1}\fi
\ifnum#1=3\ifnum#2=1\def\cycleshift{1}\fi\fi
\vc{\begin{tikzpicture}[scale=0.4]
          \outercycle{#1}{\cycleshift}
          \setcounter{Fneid}{1}         
          \foreach\i in {0,...,\vertexloopi}{   
          \pgfmathtruncatemacro{\jfrom}{\i+1}
          \foreach\j in {\jfrom,...,\vertexloop}{
          \edef\eID{\arabic{Fneid}}
          \edef\eij{\elabel[\eID]}         
            \Fne\i\j{\eij}              
	    \stepcounter{Fneid}
          }}
          \foreach\i in {0,...,\vertexloop}{\Fnv{#2}{\i}  
          }
          \end{tikzpicture}}%
\else
   #1 vertices need \expectededges{} edges but got \elabellen edges.
\fi
\fi
}

\title{Solving Large SDPs from Flag Algebra Problems}
\author{
Bernard Lidick\'{y}\thanks{Department of Mathematics, Iowa State University, Ames, IA. E-mail: {\tt lidicky@iastate.edu}. Research of this author is supported in part by NSF FRG DMS-2152490 and DMS-2554129, the Simons Foundation TSM-00013439 and Scott Hanna Professorship.}
\and Jan Volec\thanks{Department of Theoretical Computer Science, Faculty of Information Technology, Czech Technical University in Prague, Th\'akurova 9, Prague, 160 00, Czech Republic. E-mail: {\tt jan@ucw.cz}
Research of this author was supported by the grant
23-06815M of the Grant Agency of the Czech Republic.
}}
\date{\today}

\begin{document}

\maketitle
\begin{abstract}
The flag algebra method is a general tool in extremal combinatorics,
which has been useful in tackling open problems concerning graphs, hypergraphs, permutations, point sets, and other structures.
The method establishes combinatorial inequalities via sums of squares, which are usually obtained by a semidefinite relaxation of the original combinatorial problem.

A typical application of the method in graph theory starts with fixing a parameter $n$ followed by constructing a semidefinite program based on the (subgraph) structure of all $n$-vertex graphs. 
In the standard formulation before symmetry reductions, there is a coefficient constraint for every admissible $n$-vertex graph up to isomorphism.
For 3-uniform hypergraphs, there are  2,136 of them on 6 vertices and 7,013,320 on 7 vertices.
While solving the resulting program for $n=6$ is an easy task on any modern computer,
the programs for $n = 7$ seems out of reach, mainly because the SDP solvers need too much RAM. 
Yet solving such programs for as large $n$ as possible is desirable because they typically yield better bounds for the original combinatorial problems.

We bypass the aforementioned memory requirements by linear programming with cutting planes.
Although we do not guarantee that this cutting-plane procedure reaches an optimum in finitely many iterations, in many flag algebra instances it uses much less memory and yields improved bounds for well-known hypergraph Turán problems.

We use this method to obtain numerical bounds for Tur\'an densities and codegree squared densities of several 3-uniform hypergraphs including the tetrahedron, where we obtained the upper bound 0.55753 and tetrahedron without one edge with a new upper bound 0.2865.
It also allows to tackle bounds on Turán densities of 7-vertex 3-graphs using flag algebras, which was not possible before.
\end{abstract}
\def\highlight{\color{red}}

\section{Introduction}

In 2007, Razborov~\cite{Razborov} introduced the method of flag algebras for extremal combinatorics.
It allows to search for non-negative asymptotic inequalities expressible as sums of squares of small induced subgraph densities.
This is done via a semidefinite relaxation of the original (discrete) problem.
The method is quite general and can be applied to
graphs~\cite{MR4585299, MR4472294}, 
oriented graphs~\cite{MR4109639,MR4390914},
Ramsey numbers~\cite{2023_KiemPokuttaSpiegel_4colorramsey,MR4321243, lidický2024smallramseynumbersbooks},
hypergraphs~\cite{MR4608432,FALGASRAVRY2012}, permutations~\cite{MR4805530, MR3768100, Chan2020},
hypercubes~\cite{bodnár2025exactvaluesinducibilitystatistics,Balogh2014Cube}, crossing number~\cite{BaloghHill2019,Balogh2023}, SAT functions~\cite{Balogh2023SAT}, infinite trees~\cite{brosch2024gettingrootproblemsums}, 
vector spaces~\cite{2023_RueSpiegel_Radomultiplicity},
or phylogenetic trees~\cite{Alon2016} to name a few.

A typical application of flag algebras fixes an order $n$ and constructs a semidefinite relaxation using induced densities of all admissible $n$-vertex graphs or hypergraphs.
Let $f_n$ denote the number of such objects up to isomorphism. In the simplest formulation, there is one coefficient constraint for each of these $f_n$ objects, rather than one for each vertex. Larger values of $n$ typically give stronger bounds (and never worse ones).
For unrestricted $r$-uniform hypergraphs with fixed $r\geq 2$, colloquially abbreviated as \emph{$r$-graphs}, the elementary bounds
\[
 \frac{2^{\binom nr}}{n!}\leq f_n\leq 2^{\binom nr}
\]
show that $f_n=2^{\Theta(n^r)}$; in particular, this growth is superexponential in $n$. The first ten values are given in Table~\ref{tab:counts}. We note that forbidden-subgraph restrictions can substantially reduce these counts.

\begin{table}
\begin{center}
\resizebox{\textwidth}{!}{%
\begin{tabular}{c|c|c|c|c|c|c|c|c|c|c|c|c}
$r$ & OEIS & $f_0$  &  $f_1$ & $f_2$ & $f_3$ & $f_4$ & $f_5$ & $f_6$ & $f_7$ & $f_8$ & $f_9$ & $f_{10}$ \\ \hline
2          & \href{https://oeis.org/A000088}{A000088} & 1 & 1 & 2 & 4 & 11 & 34 & 156 & 1044 & 12,346 & 274,668 & 12,005,168    \\
3          & \href{https://oeis.org/A000665}{A000665}    &
1 & 1 & 1 & 2 & 5 & 34 & 2136 & 7,013,320 & $1.79 \times 10^{12}$ 
&
 $5.3 \times 10^{19}$
&
$3.66 \times 10^{29}$
\end{tabular}}
\end{center}
\caption{Number of non-isomorphic $r$-graphs up to $10$ vertices, where $r\in\{2,3\}$.}\label{tab:counts}
\end{table}

We often use the CSDP solver~\cite{Borchers1999} to solve the semidefinite programs. 
If $f_n \leq 10,000$, the resulting  semidefinite program is usually easy to solve on a personal computer.
If $f_n \leq 200,000$, it can typically be solved on a cluster in a reasonable time.
Solving instances with $f_n > 500,000$ is an exceptional challenge, although it can be done with substantial computational resources~\cite{berlinK4}.
The main limitation is the memory used by CSDP. 
In some of our work, we use MOSEK~\cite{mosek} as well, which is often significantly faster than CSDP but requires even more memory. This limits its usability for large instances.

In this paper, we describe a method for obtaining approximate SDP solutions that we applied to instances with almost $10{,}000{,}000$ coefficient constraints.
Even a partially solved relaxation at order $n$ can improve on the optimum of the relaxation at order $n-1$.
The cutting-plane method we use is already known~\cite{Krishnan2006}, even in the context of flag algebras~\cite{MR3638333}.
Nevertheless, this old method seems to work particularly well in our case.
There may be even better ways of iteratively approaching SDP solutions via LPs~\cite{Ahmadi2017}, but we have not explored these avenues.

We present our newly obtained bounds in Section~\ref{sec:results}.
Some numerical bounds were close to the conjectured exact values, so it was possible to obtain certificates for exactly ``rounded'' tight bounds.
In Section~\ref{sec:method}, we describe the method along with a toy example.  
Section~\ref{sec:fa} is dedicated to the description of flag algebras.

\section{Results on \texorpdfstring{$3$}{3}-graphs}\label{sec:results}

Our method seems to be particularly useful for Turán problems of 3-graphs. CSDP can handle calculations on six vertices, for which there are $2{,}136$ 3-graphs up to isomorphism.
However, 7 vertices are typically out of reach since there are 7,013,320 non-isomorphic 3-graphs. It seems that 6 vertices is often not quite enough to have a good approximation of the problem and stepping up to 7 vertices improves the numerical bounds significantly. 
We applied our method to unknown Tur\'an densities and codegree squared densities from the survey by Balogh, Clemen and Lidick\'y~\cite{MR4421399}
that is being occasionally updated on arXiv~\cite{survey}. 
We include the previous numerical upper bounds, and lower bounds from corresponding constructions. The previously known constructions are described in detail in the survey paper~\cite{survey}, and for the newly discovered ones, see Section~\ref{sec:constructions}.
We stopped the calculations when progress became small, without establishing optimality. Running the solver for longer may improve the reported bounds.

We test our implementation on the Tur\'an densities of the following 3-graphs. 
\begin{align*}
    K_4 &= \{123,124,134,234 \} = \binom{[4]}{3}  &  K_5 &= \binom{[5]}{3}  \\ 
    K_6 &= \binom{[6]}{3}  & K_7 &= \binom{[7]}{3} \\
    K_4^- &= \{123,124,134\}
    &
    K_5^- &= \binom{[5]}{3} \setminus\{345\}\\
    K_5^= &= \binom{[5]}{3} \setminus\{123,345\} 
    &
    K_5^{\prec} &= \binom{[5]}{3} \setminus\{234,345\} 
    \\
    C_5 &=  \{123, 234, 345 ,451, 512 \}
    &
    J_4 &= \{123, 124, 125, 134, 135, 145\}
    \\
    J_5 &=  \{123, 124, 125, 126, 134, 135, 136, 145, 146, 156\}    
    &
    F_{4^-,2} &=  \{123, 124, 134, 156, 256, 356, 456\} 
    \\
\mathbb{F} &=  \{ 123, 345, 156, 246, 147, 257, 367\}
&
    F_{4,2} &=  \{123, 124, 134, 234, 156, 256, 356, 456\}    
\end{align*}
Note that \cite{MR4421399} denotes by $F_{4,2}$ the 3-graph that we denote here by $F_{4-,2}$.\footnote{Our definition of $F_{4,2}$ is consistent with the definition of $3$-graphs $\mathcal{F}(p,q)$ in \cite[Definition 1.1]{F32Mubayi} with $p=4$ and $q=2$.}
We estimate two quantities: the Tur\'an density and its $\ell_2$ variant.

Given a 3-graph $F$, the \emph{extremal number} $\textrm{ex}(n,F)$ is the maximum number of triples in an $F$-free $n$-vertex 3-graph. The \emph{Tur\'an density} of $F$ is 
\[
\pi(F) = \lim_{n\to\infty} \frac{\textrm{ex}(n,F)}{\binom{n}{3}}.
\]
Let $G$ be an $n$-vertex 3-graph. 
The \emph{codegree} of a pair of vertices $X$ is the number of triples in $G$ containing $X$. 
Denote this number by $d_G(X)$.
We are interested in the $\ell_2$-norm of the vector of the codegrees of all pairs of vertices.
Let
\[
\textrm{co}_2(G) = \sum_{X \in \binom{V(G)}{2}} d_G(X)^2,
\]
and denote by $\textrm{exco}_2(n,F)$ the maximum of $\text{co}_2(G)$ over all $n$-vertex 3-graphs $G$ that are $F$-free. 
Finally, we define the normalized asymptotic maximum of $\text{co}_2(n,F)$, which always has a limit~\cite[Proposition 1.8]{BalCleLidmain}, as
\[
\pi_{\ell_2}(F) = \lim_{n \to \infty} \frac{\textrm{exco}_2(n,F)}{\binom{n}{2}(n-2)^2}.
\]
It seems that determining $\pi_{\ell_2}(F)$ is slightly more approachable than $\pi(F)$, since the former 
tend to have more often unique maximizers.

We summarize the bounds we obtained together with the best known bounds at the time of writing the paper in Table~\ref{tab:results}. 
For up-to-date bounds and constructions, see~\cite{survey}. 
The sizes of the corresponding semidefinite programs, as well as the resources required to solve them, are summarized in Table~\ref{tab:dats}. 


\newif\ifdecimal
\decimaltrue       

\ExplSyntaxOn
\NewDocumentCommand{\decimal}{m}
  {
    \ifdecimal
      \fp_eval:n { trunc( (#1) * 1000000 ) / 1000000 }
    \else
      #1
    \fi
  }
\ExplSyntaxOff

\ExplSyntaxOn
\NewDocumentCommand{\decimalup}{m}
  {
    \ifdecimal  
    \fp_eval:n { ( trunc( (#1) * 1000000 )+1) / 1000000 }
\else
      #1
    \fi    
  }
\ExplSyntaxOff

\newcommand{\Cfive}{
    \ifdecimal  
    0.464101
\else
    2\sqrt3-3
    \fi  
}

\begin{table}[htbp]
\centering
\resizebox{\textwidth}{!}{%
\begin{tabular}{|l||c|c|c|c|c|c|}
\hline
$F$ & $\pi(F)\geq$ & prior $\pi(F)\leq$ & our $\pi(F)\leq$ & $\pi_{\ell_2}(F)\geq$ & prior $\pi_{\ell_2}(F)\leq$ & our $\pi_{\ell_2}(F)\leq$\\
\hline\hline
$K_4^-$ & $\decimal{2/7}$~\cite{K43-extremalfrankl} & $0.28689$~\cite{flagmatic} & $0.286504$ & $\decimal{4/43}$~\cite{survey} & $0.09307$~\cite{survey} & $0.093058$\\
$K_4$ & $\decimal{5/9}$~\cite{MR177847} &
$0.55733^\star$~\cite{berlinK4}
& 
$0.557522$
& $\decimal{1/3}$~\cite{BalCleLidmain} & $\decimal{1/3}$~\cite{BalCleLidmain} & $\decimalup{(1/3)}^{\dagger}$\\
$K_5$ & $\decimal{3/4}$~\cite{Sidorenko1981SystemsOS,MR177847} & $0.7696$~\cite{flagmatic} & $0.757017$ & $\decimal{5/8}$~\cite{BalCleLidmain} & $\decimal{5/8}$~\cite{BalCleLidmain} & $\decimalup{(5/8)}^{\dagger}$\\
$K_6$ & $0.84$~\cite{Sidorenko1981SystemsOS,MR177847} & $0.8584$~\cite{flagmatic} & $0.849270$ & $0.7348$~\cite{BalCleLidmain} & $0.7536$~\cite{survey} & $0.735792$\\
$K_7$ & $\decimal{8/9}$ & $\decimal{14/15}$ & $0.897843$ & $\decimal{22/27}$~\cite{BalCleLidmain} & \pending & $0.815237$\\
$K_5^-$ & $0.602673$~\cite{HongBjarne} & $0.64209$~\cite{survey} & $0.632498$ & $0.35794$~\cite{survey} & $0.41962$~\cite{survey} & $0.406864$\\
$K_5^{=}$ & $0.58656$~\cite{survey} & $0.60962$~\cite{survey} & $0.602768$ & $0.35794$~\cite{survey} & $0.38873$~\cite{survey} & $0.382048$\\
$K_5^{\prec}$ & $\decimal{5/9}$~\cite{survey} & $0.57705$~\cite{survey} & $0.571371$ & $\decimal{1/3}$~\cite{survey} & $0.34022$~\cite{survey} & $0.333719$\\
$C_5$ & $\Cfive$~\cite{F32Mubayi} & $0.46829$~\cite{flagmatic} & $0.46503^{\ddagger}$ & $0.25194$~\cite{survey} & $0.25311$~\cite{survey} & $0.252451$\\
$J_4$ & $\decimal{1/2}$~\cite{DaisyBollobas} & $0.50409$~\cite{flagmatic} & $0.502135$ & $
\decimal{16/57}$~\cite{survey} & $0.2808$~\cite{survey} & $0.280704$\\
$J_5$ & $\decimal{9/14}$ & $0.64475$~\cite{survey} & $0.652485$ & $\decimal{27/61}$ & $0.44275$~\cite{survey} & $0.442663$\\
$F_{4^-,2}$ & $\decimal{86/175}$ & $0.493108$ & $0.493024$ & $\decimal{1552/5375}$ & $0.28888876$ & $0.288932$\\
$F_{4,2}$ &  $0.618782$  & $0.620512$ & $0.619154$ & $0.4240915$ & $0.42544754$ & $0.424094$\\
$\mathbb F$ & $\decimal{3/4}$~\cite{VeraSos} & $\decimal{3/4}$~\cite{FanoFuredi} & $\decimalup{(3/4)}^{\dagger}$ & $\decimal{5/8}$~\cite{survey} & $\decimal{5/8}$~\cite{hou2025exactturannumberfano} & $\decimalup{(5/8)}^{\dagger}$\\
\hline
\end{tabular}}
\caption{Reported bounds, with prior values and citations transcribed from the supplied local citation table where available. A dagger $\dagger$ denotes an unrounded numerical result where the exact rounding matching the lower bound is possible. A dash denotes an unreported value. The double dagger $\ddagger$ denotes the CSDP result; our method gives $0.46507$. The construction in Section~\ref{sec:constructions} gives both lower bounds for $J_5$, $F_{4^-,2}$ and $F_{4,2}$. 
$^\star$ A LEAN verified bound~\cite{jeong2026newupperboundturan} is $\pi(K_4^3) \leq 0.557808$.
    \ifdecimal
    The lower bounds typically have exact expressions coming from constructions. We opted to write numerical approximations to make it easier to see the gap between the upper and lower bounds.
    \else
Here
$2\sqrt{3}-3\approx 0.4641016$,
$22/27\approx 0.814814$,
$14/15\approx0.933333$,
$9/14\approx0.642857$, 
$27/61\approx0.442622$,
$86/175\approx0.491429$ and $1552/5375\approx0.288744$.
\fi
}
\label{tab:results}
\end{table}

\begin{table}[h!]
\begin{center}
\resizebox{\textwidth}{!}{%
\begin{tabular}{|l||c|c|c|c|c|c|c|c|c|}
\hline
 $F$ & $\ell$ &  $|\mathcal{F}_\ell| $ &  {\tt wc -l}   & dat-s & $\# B$ & Max $|B|$ &  $\sum |B|$ & our RAM  & csdp RAM\\
 \hline
 \hline
{$K_4^{-}$} & 8 &  \cellcolor{white}9,895,719    & 63,438M &  1,400 GB &202 & 5751 & 135,347 & 1,100 GB & 730,590 GB \\
{$K_4$} & 7 &  \cellcolor{white}1,295,600   &  1,427M   &   26 GB  & 52 & 990 & 15,527 & 174 GB & 12,531 GB\\
{$K_5$} & 7 &  \cellcolor{white}6,870,284   &  7,694M   &  146 GB & 72 & 990 & 34,160 & 350 GB & 351,806 GB \\
{$K_6$} & 7 &  \cellcolor{white}7,013,164   &  7,869M   &  149 GB  & 74 & 990 & 35,399 & 378 GB & 366,588 GB \\
{$K_7$} & 7 &  \cellcolor{white}7,013,319   &  7,869M   &  149 GB  & 74 & 990 & 35,400 & 280 GB& 366,604 GB\\
{$K_5^{-}$} & 7 &  \cellcolor{white}5,868,713 & 6,624M &  125 GB & 70 & 990 & 31,743 & 440 GB & 256,725 GB\\
{$K_5^{=}$} & 7 &  \cellcolor{white}4,696,659 & 5,298M &  100 GB & 68 & 990 & 29,217 & 347 GB & 164,440 GB\\
{$K_5^{\prec}$} & 7 &  \cellcolor{white}4,129,214 & 4,644M &  87 GB & 68 & 990 & 27,778 & 386 GB & 127,115 GB\\
{$C_{5}$} & 7 & 
\cellcolor{white}497,783 
&  503M  & 9.1 GB & 58 & 990 & 13,598 & 60 GB& 1,855 GB\\
{$J_4$} & 7 &  \cellcolor{white}2,982,280  & 3,368M &  63 GB  & 64 & 929 & 24,082 & 160 GB& 66,323 GB\\
{$J_5$} & 7 &  \cellcolor{white}6,761,835  & 7,606M &  144 GB  & 74 & 990 & 34,803 & 168 GB& 340,789 GB\\
{$F_{4^-,2}$} & 7 &  \cellcolor{white}2,652,406  & 2,955M &  55 GB  & 74 & 990 & 28,028 & 326 GB & 52,468 GB\\
{$F_{4,2}$} & 7 &  \cellcolor{white}5,285,756  &  5,964M  & 112 GB & 74 & 990 & 32,960 & 145 GB & 208,266 GB \\
{$\mathbb{F}$} & 7 &  \cellcolor{white}5,801,780  & 6,512M &  123 GB  & 74 & 990 & 35,400 & 131 GB & 250,904 GB\\
\hline
\end{tabular}}
\end{center}
\caption{
$|\mathcal{F}_{\ell}|$ is the number of constraints in $(P)$;
{\tt wc -l} is the number of lines in dat-s file;
dat-s is the size of the file; 
$\#B$ is the number of blocks; 
Max $|B|$ is the size of the maximum block;
$\sum |B|$ is the sum of sizes of the blocks;
RAM is the memory needed by our solver when calculating $\pi(F)$. Since the memory increases with number of iterations, it may be that smaller programs need more memory since we let them  run for longer;
CSDP RAM is the memory CSDP attempts to allocate to solve the dat-s file.
}\label{tab:dats}
\end{table}

\section{Introduction to Flag Algebras}\label{sec:fa}
The purpose of this section is to present a brief introduction to flag algebras in order to make this paper self-contained. This is not intended to be a comprehensive introduction to the method. 
The method is applicable to a variety of combinatorial structures, however, we discuss it only in the setting of graphs to make it more accessible.

For a graph $G$, denote the vertex set of $G$ by $V(G)$ and the number of vertices of $G$ by $v(G)$. For a set of vertices $X$, we denote by $G[X]$ the subgraph of $G$ induced by $X$.

Let $\sigma$ be a graph on $[k]$, where $[k]=\{1,\ldots,k\}$; in the trivial case $k=0$ it holds that $[0]=\emptyset$ and $\sigma$ is the null graph.
Such (canonically labeled) graphs $\sigma$ are called \emph{types}.
For a fixed type $\sigma$, a \emph{$\sigma$-flag} $F=(G,\theta)$ is a pair consisting of a graph $G$ and an injective map $\theta: [v(\sigma)] \to V(G)$ that is an isomorphism from $\sigma$ to $G[\operatorname{Im}(\theta)]$.
Intuitively, $F$ is a graph $G$ together with a fixed labeled copy of $\sigma$.
For a flag $F=(G,\theta)$, we write $v(F)$ to denote $v(G)$,
and, for every $\operatorname{im}(\theta) \subseteq X \subseteq V(G)$, 
we write $F[X]$ to denote the $\sigma$-flag $\left(G[X],\theta\right)$.

Two $\sigma$-flags $F_1=(G_1,\theta_1)$ and $F_2=(G_2,\theta_2)$ are isomorphic if there exists an isomorphism $\varphi$  from $V(G_1)$ to $V(G_2)$ such that for all  $i \in V(\sigma): \varphi(\theta_1(i)) = \theta_2(i)$.
In other words, $\varphi$ is an isomorphism from $F_1$ to $F_2$ that maps the labeled copy of $\sigma$ in $F_1$ to the labeled copy of $\sigma$ in $F_2$, and the induced automorphism of $\sigma$ is the identity.
Denote the family of all non-isomorphic $\sigma$-flags with $n$ vertices by $\mathcal{F}^\sigma_n$, and let $\mathcal{F}^\sigma := \bigcup_{n\in\mathbb{N}}\mathcal{F}^\sigma_n$.

Given two $\sigma$-flags $F_1=(G_1,\theta_1)$ and $F_2=(G_2,\theta_2)$, let
$P(F_1,F_2)$ be the probability that $F_2[X]$ is a $\sigma$-flag isomorphic to $F_1$, where
$X$ is a uniformly random $v(F_1)$-vertex subset of $V(F_2)$ such that $\operatorname{Im}(\theta_2) \subseteq X$.
Clearly, $P(F_1,F_2) \in [0,1]$, and also note that $P(F_1,F_2)=0$ in case $v(F_1) > v(F_2)$.

For three $\sigma$-flags 
$F_1=(G_1,\theta_1)$,  
$F_2=(G_2,\theta_2)$
and
$F_3=(G_3,\theta_3)$ such that $v(F_3) = v(F_1) + v(F_2)-v(\sigma)$,
we define the \emph{pair density $P(F_1,F_2; F_3)$} to be
the probability that $F_1 \cong F_3[V(G_3) \setminus X]$
and
$F_2 \cong F_3[X \cup \operatorname{Im}(\theta_3)]$,
where $X$ is a uniformly random subset of $V(G_3)\setminus \operatorname{Im}(\theta_3)$ of size $v(F_2)-v(\sigma)$.
Note that $v(F_3) - |X| = v(F_1)$ and, analogously to the previous paragraph, it holds that $P(F_1,F_2; F_3) \in [0,1]$. 

Given a type $\sigma$, let $\mathbb{R}\mathcal{F}^\sigma$ be the set of all formal finite linear combinations of $\sigma$-flags with real coefficients.
The addition on $\mathbb{R}\mathcal{F}^\sigma$ as well as the multiplication of the elements of $\mathbb{R}\mathcal{F}^\sigma$ by a real number are defined naturally.

Consider the set $\mathcal{K^\sigma}$ that is equal to the span of all linear combinations of the form
\[
F - \sum_{H \in \mathcal{F}_{v(F)+1}^\sigma} P(F,H) \cdot H\,, \quad  F\in \mathcal{F}^\sigma,
\]
and observe that an easy double-counting argument reveals that for two $\sigma$-flags $F$ and $G$ with $v(F)<v(G)$, it holds that 
\[
P(F,G) - \sum_{H \in \mathcal{F}_{v(F)+1}^\sigma} P(F,H) \cdot P(H,G) = 0.
\]

Let $\mathcal{A}^\sigma=\mathbb{R}\mathcal{F}^\sigma/\mathcal{K^\sigma}$ be the corresponding quotient vector space. 
To finish the definition of the algebra $\mathcal{A}^\sigma$, we linearly extend the following definition of the product of any two $\sigma$-flags $F_1$ and $F_2$:
\[
F_1 \times F_2 = \sum_{F_3 \in \mathcal{F}_{v(F_1)+v(F_2)-v(\sigma)}^\sigma} P(F_1,F_2;F_3) \cdot F_3.
\]
Razborov proved that the extended product $\times$ on $\mathbb{R}\mathcal{F}^\sigma$ is compatible with the factorization, and thus $\mathcal{A}^\sigma$ is indeed an algebra.
For brevity, in case $\sigma$ is the null graph, we identify the $\sigma$-flag with the underlying graph, and we also write $\mathcal{A}$ instead of $\mathcal{A}^\sigma$.

Given a $\sigma$-flag $F=(G,\theta)$, let $q_F$ be the probability that a uniformly random injective map $\psi: [v(\sigma)] \to V(G)$ satisfies $F \cong (G, \psi)$.
As in Razborov, for every type $\sigma$, we define the operator $\llbracket . \rrbracket_\sigma : \mathcal{A}^\sigma \to \mathcal{A}$ as the unique linear extension of the following map from $\mathcal{F}^\sigma$ to $\mathcal{F}$:
\[
\llbracket F \rrbracket_\sigma = q_F \cdot G, \quad \mbox{where $F = (G,\theta)$}.
\]
For brevity, we drop $\sigma$ in the subscript of $\llbracket \cdot \rrbracket$ whenever the type is clear from the context.

The algebras connect to graphs using convergent sequences and homomorphisms. 
A sequence of graphs $(G_n)$ of increasing orders is called \emph{convergent} if for every finite graph $H$
\[
\lim_{n\to\infty} P(H,G_n)
\]
exists.
Observe that by compactness, every sequence of graphs has a convergent subsequence.
The limit object of a convergent sequence can be naturally represented by a map $\phi:\mathcal F\to[0,1]$. Razborov proved that the map extends to a unital algebra homomorphism $\phi:\mathcal A\to\mathbb R$ that is \emph{positive}, i.e., it satisfies $\phi(F)\ge0$ for every $F \in \mathcal{F}$. On the other hand, a sampling argument reveals that every such homomorphism arises from at least one convergent sequence of graphs of increasing orders.

Let $\operatorname{Hom}^+(\mathcal A,\mathbb R)$ be the set of all unital algebra homomorphisms from $\mathcal{A}$ to $\mathbb{R}$.
A key feature of the flag algebra framework is that for every
 $\phi \in \operatorname{Hom}^+(\mathcal{A},\mathbb{R})$
 and for every $f \in \mathcal{A}^\sigma$, it holds that
 \[
 \phi \left(  \left\llbracket  f^2 \right\rrbracket  \right) \geq 0.
 \]
This allows us to generate many subgraph-counting inequalities that are valid for all convergent sequences of graphs simply using squares of linear combinations of $\sigma$-flags.

\subsection{The Turán densities $\pi(\cdot)$ and $\pi_{\ell_2}(\cdot)$ in Flag Algebras}

Given a question on the Turán density of a fixed hypergraph $F$,
it is straightforward to design a hiearchy of its flag-algebraic semidefinite relaxations.
Indeed, given an integer $k$, consider all inequalities of the form $\left\llbracket f^2 \right\rrbracket$ that can be rewritten as a linear combination of $k$-vertex subgraph densities.

A moment of thought reveals that bounding $\pi_{\ell_2}(F)$ is as easy task as bounding $\pi(F)$, since the objective can be expressed as a linear combination of densities of 4-vertex 3-graphs.
Indeed, consider any convergent sequence of 3-graphs $(G_n)$ whose orders $N_n$ tend to infinity, and let $\phi$ be the corresponding limit homomorphism.
The scaled limit of $\textup{co}_2(G_n)$ can be rewritten using $\phi$ as follows:

\begin{align*}
\lim_{n \to \infty} \frac{\textup{co}_2(G_n)}{\binom{N_n}{2}(N_n-2)^2}
=
\phi\left(
\left\llbracket
\left(
\vc{
\begin{tikzpicture}
\draw
(0,0) coordinate(1) node[label=below:1](a){}
(1,0) coordinate(2) node[label=below:2](b){}
(0.5,1) coordinate(3) node[vtx](c){}
;
\draw[hyperedge] (1) to[out=40,in=130] (2) to[out=120,in=280] (3) to[out=260,in=50] (1);
\draw
(1) node[labeled_vertex]{}
(2) node[labeled_vertex]{}
;
\end{tikzpicture}
}
\right)^2
\right\rrbracket
\right)
=
\phi\left(
\frac{1}{6}
\vc{
\begin{tikzpicture}
\draw
(0,0) coordinate(1) node[vtx](a){}
(1,0) coordinate(2) node[vtx](b){}
(1,1) coordinate(3) node[vtx](c){}
(0,1) coordinate(4) node[vtx](d){}
;
\draw[hyperedge] (1) to[out=30,in=90] (2) to[out=100,in=260] (3) to[out=260,in=30] (1);
\draw[hyperedge] (1) to[out=90,in=150] (2) to[out=150,in=280] (4) to[out=280,in=90] (1);
\end{tikzpicture}
}
+
\frac{1}{2}
\vc{
\begin{tikzpicture}[scale=0.55]
\foreach \i in {1,2,3,4}{
\draw (90*\i-45:1.5) coordinate(\i);
}
\foreach \r in {0-45,90-45,180-45}{
\begin{scope}[rotate=\r]
\draw[hyperedge] 
(0:1.5) to[out=140,in=275,looseness=1.2] (90:1.5)  to[out=265,in=40,looseness=1.2] (180:1.5)  to[out=35,in=145,looseness=1.2] (0:1.5) 
;
\end{scope}
}
\draw 
(1) node[vtx,label=right:{\tiny }]{}
(2) node[vtx,label=left:{\tiny }]{}
(3) node[vtx,label=left:{\tiny }]{}
(4) node[vtx,label=right:{\tiny }]{}
;
\end{tikzpicture}
}
+
\vc{
\begin{tikzpicture}[scale=0.55]
\foreach \i in {1,2,3,4}{
\draw (90*\i-45:1.5) coordinate(\i);
}
\foreach \r in {0-45,90-45,180-45,270-45}{
\begin{scope}[rotate=\r]
\draw[hyperedge] 
(0:1.5) to[out=140,in=275,looseness=1.2] (90:1.5)  to[out=265,in=40,looseness=1.2] (180:1.5)  to[out=35,in=145,looseness=1.2] (0:1.5) 
;
\end{scope}
}
\draw 
(1) node[vtx,label=right:{\tiny }]{}
(2) node[vtx,label=left:{\tiny }]{}
(3) node[vtx,label=left:{\tiny }]{}
(4) node[vtx,label=right:{\tiny }]{}
;
\end{tikzpicture}
}
\right).
\end{align*}

\section{Solving large Flag Algebra relaxations}\label{sec:method}

The method we describe below is known in the extremal combinatorics community. 
It was used in the early days of flag algebras by Hladk\'y, Kr\'al' and Norin~\cite{MR3638333} before the community switched to semidefinite programming solvers such as CSDP~\cite{Borchers1999} and SDPA~\cite{SDPA1, SDPA2, SDPA-GMP, SDPA3}. 
When looking for references, we have also found that the method we use is described on a GUROBI~\cite{gurobi} support webpage~\cite{GurobiCommunity}.
A similar method is occasionally used by LLMs when producing flag algebra programs, and, in particular, it was used recently by Jeong, Park, Im, Lee, and Yang~\cite{jeong2026newupperboundturan} on Tur\'an's tetrahedron problem.

We use the following equality-form of a semidefinite program and its dual. All matrices are real and symmetric, and $A\cdot X=\operatorname{tr}(A^TX)$ is the trace inner product. The constraints are indexed by $i=1,\ldots,m$. 

\begin{align}
(P)
\begin{cases}
\text{max}  & C \cdot X \\
\text{s.t}     &  A_i \cdot X = b_i \\
& X \succcurlyeq 0
\end{cases}
&&
(D)
\begin{cases}
\text{min}  & b^T y \\
\text{s.t}     &  \sum_i y_i A_i  - C  \succcurlyeq 0 \\
& y \in \mathbb{R}^m
\end{cases}
\end{align}

Typically, $X \succcurlyeq 0$ is a block diagonal matrix. 
We attempt to solve $(D)$ by forgetting all blocks that are not \emph{trivial}, i.e., of size $1 \times 1$.
Let $A_i'$ and $C'$ be the submatrices of $A_i$ and $C$, respectively, whose entries correspond only to trivial blocks in $X$, i.e., $A_i'$ and $C'$ are diagonal matrices.

\begin{align}
(D')
\begin{cases}
\text{min}  & b^T y \\
\text{s.t}     & \sum_i y_i A'_i  - C'  \geq 0 \\
& y \in \mathbb{R}^m
\end{cases}
\end{align}

The relaxation $(D')$ is a linear program. Assuming that it has a finite attained optimum, let $y^\star$ be an optimal solution. If
\[
M(y^\star):=\sum_i y_i^\star A_i-C\succeq0,
\]
then $y^\star$ is feasible and optimal for $(D)$, since $(D')$ is its relaxation. Otherwise, for each block $B(y^\star)$ that is not positive-semidefinite, choose a unit eigenvector $v$ for its smallest (negative) eigenvalue, and add to $(D')$ the linear inequality
\[
v^T B(y)v\geq0.
\]
Here $v$ is fixed while $y$ remains a variable. Equivalently, extend $v$ by zeros outside that block and impose $v^TM(y)v\geq0$.
Resolve the new linear program and repeat. Each cut excludes the current solution while preserving every feasible solution of $(D)$.
The relaxation's minimum is nondecreasing, but neither strict objective improvement nor finite termination are guaranteed. The initial relaxation can also be unbounded for a general SDP; in our density applications, however, the normalization and the nonnegativity constraints yield a compact feasible region.

\subsection{Application to flag algebras}

In the application to flag algebras,  $y$ in $(D)$ can be thought as subgraph densities $\phi(F)$ for $F \in \mathcal{F}_k$ of some $\phi \in \operatorname{Hom}^+(\mathcal{A},\mathbb{R})$.
This naturally gives linear constraints 
$\sum_i y_i = 1$ and $y_i \geq 0$ for all $i \in [f_k]$. 
We use the density maximization formulation below, whose objective has the opposite direction from the general minimization formulation $(D)$. The nontrivial matrix blocks are homogeneous in $y$, i.e., in the notation above, their corresponding blocks in $C$ are all zero. Equality constraints such as the normalization are retained explicitly. 

Each non-trivial block corresponds to a type $\sigma$ and a vector of labeled flags $F^\sigma$. 
If $\sum_i y_i A_i^\sigma  - C^\sigma  \succcurlyeq 0 $ does not hold, then we have a negative eigenvalue in the block. Choose an eigenvector $v$ corresponding to the smallest eigenvalue. 
Consider the flag algebra inequality
\[
0 \leq \phi\!\left(\llbracket (v^T F^\sigma)^2 \rrbracket\right) = \sum_i c_i\phi(F_i),
\]
which holds for every $\phi\in \operatorname{Hom}^+(\mathcal{A},\mathbb{R})$,
and corresponds to the following additional linear constraint to $(D)$:
\[
0 \leq \sum c_i y_i. 
\]

We now illustrate the method by establishing the flag algebra form of Mantel's theorem.

\begin{theorem}[Mantel~\cite{zbMATH02643703}]
    Every triangle-free graph on $N$ vertices has at most $\lfloor N^2/4\rfloor$ edges. Consequently, the density of triangle-free graphs is at most $1/2+o(1)$.
\end{theorem}

We calculate with $3$-vertex subgraph densities, thus the graphs we consider are
$\mathcal{F}_3 =\left(\Fuuu111, \Fuuu112,\Fuuu221\right)$. 
We work in the triangle-free flag algebra, so any flag containing a triangle is omitted.
Let $\phi$ be any limit homomorphism of a (convergent) sequence of triangle-free graphs.
Clearly,
\begin{align}
\phi\left(\Fuu2\right)
=
\phi\left(0\Fuuu111+\frac{1}{3}\Fuuu112+\frac{2}{3}\Fuuu221\right).
\label{eq:a}
\end{align}

Consider the flag vector $F = \left(\Flu2,\Flu1\right)^T$ and a positive semidefinite matrix with entries $a,b,c$. In what follows, writing $f\geq0$ for an algebra element $f \in \mathcal{A}$ means that $\phi(f)\geq0$ for every admissible homomorphism $\phi$.
\begin{align}
0 \leq \left\llbracket
\left(\Flu2,\Flu1\right)
\begin{pmatrix}
         a &           b   \\
         b &           c 
\end{pmatrix}
\begin{pmatrix}
         \Flu2  \\
        \Flu1 
\end{pmatrix}
\right\rrbracket=
\frac{3c}{3}\Fuuu111 + \frac{c+2b}{3}\Fuuu112 + \frac{a+2b}{3}\Fuuu221
\label{eq:b}
\end{align}

Adding the nonnegative evaluation of \eqref{eq:b} to \eqref{eq:a} and multiplying by $3$ gives
\begin{align*}
3\phi\left(\Fuu2\right)
&\leq\phi\left(3c\Fuuu111+(1+c+2b)\Fuuu112+(2+a+2b)\Fuuu221\right)\\
&\leq\max\{3c,1+c+2b,2+a+2b\},
\end{align*}
since the three $3$-vertex induced subgraph densities are nonnegative and sum to one.

The best upper bound obtainable in this way can be written as a semidefinite program.

\begin{align}
(P)
\begin{cases}
\text{min}  & t  \\
\text{s.t}     & 
 \underbrace{\begin{pmatrix}
         0 &           0  \\ 
         0 &           3 
\end{pmatrix}}_{=A_1}
  \cdot 
   \underbrace{\begin{pmatrix}
         a &           b   \\
         b &           c 
\end{pmatrix}}_{=X}
+
0
   \leq t \\
 & 
 \underbrace{\begin{pmatrix}
         0 &           1  \\ 
         1 &           1 
\end{pmatrix}}_{=A_2}
  \cdot 
\begin{pmatrix}
         a &           b   \\
         b &           c 
\end{pmatrix}
+
1
   \leq t \\ 
   &
 \underbrace{\begin{pmatrix}
         1 &           1  \\ 
         1 &           0 
\end{pmatrix}}_{=A_3}
  \cdot 
\begin{pmatrix}
         a &           b   \\
         b &           c \\
\end{pmatrix}
+
2
   \leq t \\        
& \begin{pmatrix}
         a &           b   \\
         b &           c \\
\end{pmatrix} \succcurlyeq 0
\end{cases}
&&
(D)
\begin{cases}
\text{max}  & 0 y_1 + y_2 + 2y_3 \\
\text{s.t}     &  \sum_i y_i A_i    \succcurlyeq 0 \\
& y_1+y_2+y_3 = 1\\
& 0 \leq y_i 
\end{cases}
\end{align}

Here $t$ is unrestricted, so its dual constraint is $y_1+y_2+y_3=1$; the primal inequalities give $y_i\geq0$. The objective is three times the edge density.

First we ignore the positive semidefinite constraint and solve 
\begin{align}
(D_1)
\begin{cases}
\text{max}  & 0 y_1 + y_2 + 2y_3 \\
\text{s.t}     & y_1+y_2+y_3 = 1\\
& 0 \leq y_i
\end{cases}
\end{align}
The obvious solution is $(0,0,1)$.
The positive semidefinite constraint claims
$
A := 0A_1 + 0A_2+A_3 \succeq 0.
$
However,
\[
 A = A_3 =  \begin{pmatrix}
         1 &           1  \\ 
         1 &           0 
\end{pmatrix}
\]
has numerical eigenvalues $1.618$ and $-0.618$, where the second one is negative with eigenvector $(1-\sqrt{17},4)$, or, numerically, with $v=(-0.525, 0.850)$.

Since $v$ is (numerically) a negative eigenvector,
\begin{align*}
0 > v^TAv =  \left( \sum_i y_i A_i\right) v  = 
\sum_i y_i \left(v^TA_i v\right).
\end{align*}
Therefore, adding
\[ \sum_i y_i \left(v^TA_i v\right) \geq 0 \]
to $(D_1)$ yields another linear program.
In the flag algebra language, it corresponds to the following inequality:
\begin{align*}
0 \leq \left\llbracket \left(  v_1 \Flu2 + v_2 \Flu1  \right)^2\right\rrbracket
&=
 \left\llbracket 
v_1^2\Fluu221+v_1v_2\Fluu121+ v_1v_2\Fluu122+ v_2^2\Fluu111
+ v_2^2\Fluu112
\right\rrbracket\\
&=
v_2^2\Fuuu111 
+  
\frac{2v_1v_2+v_2^2}{3}\Fuuu112
+
\frac{v_1^2+2v_1v_2}{3}\Fuuu221.
\end{align*}
Plugging in the values of $v$ yields a new constraint
$
2.171y_1 - 0.171y_2 -0.618y_3 \geq 0
$,
and in turn we obtain the following new linear program:
\begin{align}
(D_2)
\begin{cases}
\text{max}  & 0 y_1 + y_2 + 2y_3 \\
\text{s.t}         & y_1+y_2+y_3 = 1\\
& 2.171y_1 - 0.171y_2 -0.618y_3 \geq 0 \\
& 0 \leq y_i 
\end{cases}
\end{align}
This time, the optimal solution of ($D_2)$ is 
$y_1 = 0.221,
y_2 = 0.00,
y_3 = 0.778
$
with value $1.556$, and
\[
A_2= 
 \begin{pmatrix}
         0.778 &        0.778 \\ 
          0.778 &        0.664
 \end{pmatrix}
\]

The eigenvalue $-0.058$ of $A_2$ has an eigenvector $(-0.681, 0.732)$, hence we establish a new (useful) constraint
$1.609y_1  -0.4609y_2 -0.533 y_3 \geq 0$.

After adding the last constraint, the objective drops to $1.501$.
Applying the method one more time yields a constraint $ 1.503y_1  -0.498y_2 -0.501y_3 \geq 0$
(via an eigenvector $(-0.7062,0.70797)$).
The resulting program $(D_4)$ is below.
\begin{align*}
(D_4)
\begin{cases}
\text{max}  & 0 y_1 + y_2 + 2y_3 \\
\text{s.t}     & y_1+y_2+y_3 = 1\\
                  & 2.171y_1 - 0.171y_2 -0.618y_3 \geq 0 \\
                  &1.609y_1  -0.4609y_2 -0.533 y_3 \geq 0\\
		 &1.503y_1  -0.498y_2    -0.501y_3 \geq 0 \\
& 0 \leq y_i 
\end{cases}
\end{align*}
The solution of $(D_4)$ is 1.500002283054241 with 
$y_1= 0.24999885847287948$,
$y_2= 0.0$
$y_3= 0.7500011415271205$
and eigenvalues are
$1.5000000000034748$ and $-0.000002283057$.
It is quite close to the optimum value $1.5$ (with the densities $y_1=1/4$, $y_2=0$ and $y_3=3/4$,
which, indeed, corresponds to the combinatorial solution in Mantel's theorem).

\subsection{Certificate of correctness}\label{sec:certificates}

The method described in this section generates simple certificates of the established bounds by inspecting the duals of the linear programs. 
Here we show and example for the solution of $(D_4)$. 
In the dual solution of $(D_4)$, out of the three variables corresponding to the inequalities in $(D_4)$, only the last one is positive. Specifically, it is 0.997536.
This produces the following proof written in flag algebras as
\begin{align*}
3 \phi\left(\Fuu2\right) &\leq 3 \phi\left(\Fuu2\right) +  \phi\left(0.997536 \cdot 3\left\llbracket \left( -0.7062\Flu1 + 0.70797 \Flu2 \right)^2 \right \rrbracket\right)
\\
&=
\phi\left(
1.4999  \Fuuu111  
+0.5025  \Fuuu112  
+1.50001 \Fuuu122 
\right)
\leq  1.50001.
\end{align*}
Notice we obtained an upper bound $1.50001$, which is not equal $1.5$, to the optimal value of $(D)$. 

In general, the certificate from our solution procedure can be written as a linear combination of sum of squares, where the coefficients in the squares come from eigenvectors corresponding to negative eigenvalues and the coefficients in the linear combination come from the dual of the linear programs.
In particular, the certificate itself does not need the semidefinite program $(P)$, which makes the certificate quite compact.

\subsection{Implementation}

We have two independent implementations of the above method in C and C++ (each author has implemented one of them). 
We also have a mock implementation in Python, however, as we apply the method to large semidefinite programs, it was useful to implement the actual solver in C/C++ in order to reduce the memory footprint.
Solutions to all the linear programs were obtained using MOSEK or GUROBI. 
Our implementation supports both solvers.

The method seems to perform very well on instances where only a few coefficients in the objective function are large.
This seems to often be the case for Tur\'an type problems. 

Unlike CSDP, the memory needed during the calculation and the time per iteration increase as more constraints are generated and added to the linear program.
Therefore, it would be useful to keep only constraints that are needed.
In our experiments, not all linear constraints are typically needed when solving the linear program. 
However, experiments showed that a constraint that is not needed in one iteration may be useful in a later iteration.
Our implementation adds constraints during a run and saves a checkpoint containing the constraints retained from the last iteration. 
When a run has stopped because of a time limit or insufficient memory, we have restarted from the smaller checkpoint.
The implementation in C++ also produces certificates described in Section~\ref{sec:certificates}.

Our implementation of the solver is available at \url{https://github.com/BLL4/lpsdp},
and the certificates are available at
\url{https://lidicky.name/pub/lpsdp/} or 
\url{https://honza.ucw.cz/proj/lpsdp/}.

\subsection{Software for verification}

In order to verify the claimed flag algebra upper bounds, it is enough to evaluate respective feasible dual solutions, which in this case are stored as certificates described in Section~\ref{sec:certificates}.
To do so, we have implemented a C program that generates all the flag algebra coefficients in parallel using OpenMP~\cite{OpenMP}, and then it uses MPFI~\cite{MPFI} and MPFR~\cite{MPFR} libraries for floating-point interval arithmetic in order to compute the corresponding upper bounds on $\pi(\cdot)$ and $\pi_{\ell_2}(\cdot)$.

\section{Constructions}\label{sec:constructions}
\begin{figure}[h]
\begin{center}
  \begin{tabular}{ | l|| c |  c | c |   } 

\hline
$F_{4^-,2}$ &
\vc{ 
\begin{tikzpicture}[flag_pic]\outercycleH{6}{0}
\drawHEthree{6}{0}{0}{1}{2}
\drawHEthree{6}{1}{0}{1}{3}
\drawHEthree{6}{2}{0}{2}{3}
\drawHEthree{6}{3}{0}{4}{5}
\drawHEthree{6}{4}{1}{4}{5}
\drawHEthree{6}{5}{2}{4}{5}
\drawHEthree{6}{6}{3}{4}{5}
\end{tikzpicture}
}
&
\vc{
\begin{tikzpicture}
\draw
(0,-0.7) coordinate(1) node[vtx,label=below:{\tiny $1$}](a){}
(0,-0.2) coordinate(2) node[vtx,label=below:{\tiny $2$}](b){}
(0, 0.2) coordinate(6) node[vtx,label=above:{\tiny $3$}](b){}
(0,0.7) coordinate(3) node[vtx,label=left:{\tiny $4$}](c){}
(1,-0.5) coordinate(4) node[vtx,label=right:{\tiny $6$}](d){}
(1,0.5) coordinate(5) node[vtx,label=right:{\tiny $5$}](d){}
;
\draw[hyperedge] (1) to[bend left] (2) to[bend left] (3) to[bend right] (1);
\draw[hyperedge] (1) to[bend left=40] (6) to[bend left=40] (3) to[bend right=40] (1);
\draw[hyperedge] (1) to[bend left=70,looseness=2] (2) to[bend left=70,looseness=2] (6) to[bend right=70,looseness=2] (1);
\draw[hyperedge] (1) to[out=0,in=190] (4) to[out=190,in=250] (5) to[out=250,in=0] (1);
\draw[hyperedge] (2) to[out=0,in=140] (4) to[out=140,in=210,looseness=1.5] (5) to[out=210,in=0] (2);
\draw[hyperedge] (6) to[out=0,in=140] (4) to[out=140,in=210,looseness=1.5] (5) to[out=210,in=0] (6);
\draw[hyperedge] (3) to[out=0,in=120,looseness=1.2] (4) to[out=120,in=170,looseness=1.2] (5) to[out=170,in=0,looseness=1.2] (3);
\end{tikzpicture}
} 
\\\hline
$F_{4,2}$ &
\vc{ 
  \begin{tikzpicture}[flag_pic]
  \outercycleH{6}{0}
\drawHEthree{6}{0}{0}{1}{2}
\drawHEthree{6}{1}{0}{1}{3}
\drawHEthree{6}{2}{0}{2}{3}
\drawHEthree{6}{3}{1}{2}{3}
\drawHEthree{6}{4}{0}{4}{5}
\drawHEthree{6}{5}{1}{4}{5}
\drawHEthree{6}{6}{2}{4}{5}
\drawHEthree{6}{7}{3}{4}{5}
\end{tikzpicture} } 
&
\vc{
\begin{tikzpicture}
\draw
(0,-0.7) coordinate(1) node[vtx,label=below:{\tiny $1$}](a){}
(0,-0.2) coordinate(2) node[vtx,label=below:{\tiny $2$}](b){}
(0, 0.2) coordinate(6) node[vtx,label=above:{\tiny $3$}](b){}
(0,0.7) coordinate(3) node[vtx,label=left:{\tiny $4$}](c){}
(1,-0.5) coordinate(4) node[vtx,label=right:{\tiny $6$}](d){}
(1,0.5) coordinate(5) node[vtx,label=right:{\tiny $5$}](d){}
;
\draw[hyperedge] (1) to[bend left] (2) to[bend left] (3) to[bend right] (1);
\draw[hyperedge] (1) to[bend left=40] (6) to[bend left=40] (3) to[bend right=40] (1);
\draw[hyperedge] (1) to[bend left=70,looseness=2] (2) to[bend left=70,looseness=2] (6) to[bend right=70,looseness=2] (1);
\draw[hyperedge] (2) to[bend left=70,looseness=2] (6) to[bend left=70,looseness=2] (3) to[bend right=70,looseness=2] (2);
\draw[hyperedge] (1) to[out=0,in=190] (4) to[out=190,in=250] (5) to[out=250,in=0] (1);
\draw[hyperedge] (2) to[out=0,in=140] (4) to[out=140,in=210,looseness=1.5] (5) to[out=210,in=0] (2);
\draw[hyperedge] (6) to[out=0,in=140] (4) to[out=140,in=210,looseness=1.5] (5) to[out=210,in=0] (6);
\draw[hyperedge] (3) to[out=0,in=120,looseness=1.2] (4) to[out=120,in=170,looseness=1.2] (5) to[out=170,in=0,looseness=1.2] (3);
\end{tikzpicture}
} 
\\   \hline
   
   $J_5$
   &
\vc{\begin{tikzpicture}\outercycleH{7}{6}
\drawhyperedge{0}{6}
\drawhypervertex{0}{0}
\drawhypervertex{1}{0}
\drawhypervertex{2}{0}
\drawhyperedge{1}{6}
\drawhypervertex{0}{1}
\drawhypervertex{1}{1}
\drawhypervertex{3}{1}
\drawhyperedge{2}{6}
\drawhypervertex{0}{2}
\drawhypervertex{1}{2}
\drawhypervertex{4}{2}
\drawhyperedge{3}{6}
\drawhypervertex{0}{3}
\drawhypervertex{1}{3}
\drawhypervertex{5}{3}
\drawhyperedge{4}{6}
\drawhypervertex{0}{4}
\drawhypervertex{2}{4}
\drawhypervertex{3}{4}
\drawhyperedge{5}{6}
\drawhypervertex{0}{5}
\drawhypervertex{2}{5}
\drawhypervertex{4}{5}
\drawhyperedge{6}{6}
\drawhypervertex{0}{6}
\drawhypervertex{2}{6}
\drawhypervertex{5}{6}
\drawhyperedge{7}{6}
\drawhypervertex{0}{7}
\drawhypervertex{3}{7}
\drawhypervertex{4}{7}
\drawhyperedge{8}{6}
\drawhypervertex{0}{8}
\drawhypervertex{3}{8}
\drawhypervertex{5}{8}
\drawhyperedge{9}{6}
\drawhypervertex{0}{9}
\drawhypervertex{4}{9}
\drawhypervertex{5}{9}
\end{tikzpicture} 
}
&
\vc{
\begin{tikzpicture}
\draw
(-2,0) coordinate(1) node[vtx,label=left:{\tiny $1$}](1){}
\foreach \x in {2,...,6}
{
(-144+72*\x:1) coordinate(\x) node[vtx](v\x){}
}
;
\draw[hyperedge,opacity=0.2] (1)--(2)--(3)--(1);
\draw[hyperedge,opacity=0.2] (1)--(2)--(4)--(1);
\draw[hyperedge,opacity=0.2] (1)--(2)--(5)--(1);
\draw[hyperedge,opacity=0.2] (1)--(2)--(6)--(1);
\draw[hyperedge,opacity=0.2] (1)--(3)--(4)--(1);
\draw[hyperedge,opacity=0.2] (1)--(3)--(5)--(1);
\draw[hyperedge,opacity=0.2] (1)--(3)--(6)--(1);
\draw[hyperedge,opacity=0.2] (1)--(4)--(5)--(1);
\draw[hyperedge,opacity=0.2] (1)--(4)--(6)--(1);
\draw[hyperedge,opacity=0.2] (1)--(5)--(6)--(1);
\draw [line width = 0.5pt] (2)--(3)--(4)--(5)--(6)--(2)--(4)--(6)--(3)--(5)--(2)  ;
\draw (-2,0) coordinate(1) node[vtx](1){}
\foreach \x in {2,...,6}
{
(-144+72*\x:1) coordinate(\x) node[vtx](v\x){}
};
\end{tikzpicture}
}
   \\
   \hline
\end{tabular}
\end{center}
\caption{3-graphs $F_{4^-,2}$, $F_{4,2}$ and $J_5$.}\label{fig-3graphs}
\end{figure}
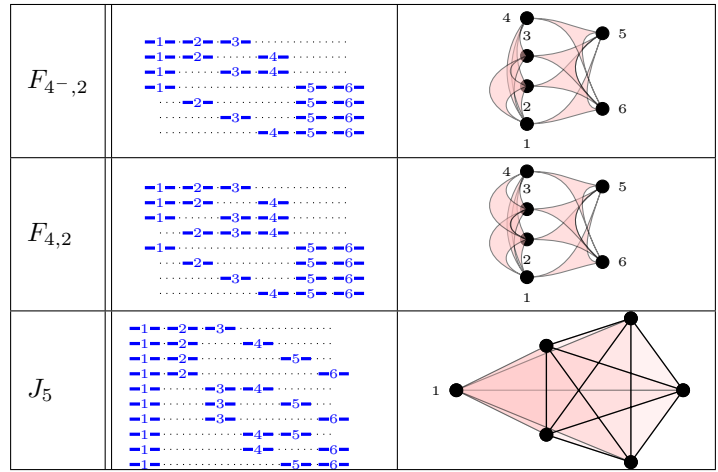

This section presents three respective constructions of $F$-free $3$-graphs for $F \in \left\{  F_{4^-,2}, F_{4,2},J_5\right\}$. Their visualizations are in Figure~\ref{fig-3graphs}.
These give explicit lower bounds for $\pi(F_{4^-,2})$, $\pi_{\ell_2}(F_{4^-,2})$, 
$\pi(F_{4,2})$, $\pi_{\ell_2}(F_{4,2})$,
$\pi(J_{5})$ and  $\pi_{\ell_2}(J_{5})$.

All three constructions were discovered by ChatGPT-Astra.

\subsection{A construction for $\pi\left(F_{4-,2}\right)\ge \frac{86}{175}$ and $\ell_2\left(F_{4-,2}\right) \ge \frac{1552}{5375}$}
\label{sec:constructions:F4-2}

\begin{figure}[ht]
\begin{center}
\tikzset{
  bigvtx/.style={fill=white,draw},
  hyperedgeI/.style={fill=green,opacity=0.5,draw=black},      
  hyperedgeA/.style={fill=red!40,opacity=0.35,draw=black},  
  hyperedgeB/.style={fill=blue!35,opacity=0.35,draw=black}, 
}
\newcommand\hv[3]{\path (c#1) ++({270+72*#1+#2}:\dout) coordinate (#3);}

\vc{
\begin{tikzpicture}
\def\R{2}       
\def\rc{0.8}    
\def\din{0.1}  
\def\rin{0.64}  
\def\dout{0.5}  
\def\bendA{10}  
\def\bendB{8}   
\def\sAfirst{-62.5} \def\sAlast{62.5}  \def\sAmid{-12.5}
\def\sBapex{12.5} \def\sBfirst{-37.5}  \def\sBsecond{37.5} 

\foreach \i in {0,...,4}
  \draw[bigvtx] ({90+72*\i}:\R) coordinate (c\i) circle (\rc);

\foreach \i in {0,...,4}{
  \path (c\i) ++({90+72*\i}:\din) coordinate (m)
        (m) +({90+72*\i}:\rin)     coordinate (1)
            +({90+72*\i+120}:\rin) coordinate (2)
            +({90+72*\i+240}:\rin) coordinate (3);
  \draw[hyperedgeI] (1) to[bend left] (2) to[bend left] (3) to[bend left] (1);
}

\foreach \k in {0,...,4}{
  \pgfmathtruncatemacro\ka{mod(\k+1,5)}
  \pgfmathtruncatemacro\kb{mod(\k+2,5)}
  \pgfmathtruncatemacro\kc{mod(\k+3,5)}
  \hv{\k}{\sAfirst}{p1} \hv{\ka}{\sAmid}{p2} \hv{\kb}{\sAlast}{p3}
  \draw[hyperedgeA] (p1) to[bend left=\bendA] (p2)
                         to[bend left=\bendA] (p3)
                         to[bend left=\bendA] (p1);
  \hv{\k}{\sBfirst}{p1} \hv{\ka}{\sBsecond}{p2} \hv{\kc}{\sBapex}{p3}
  \draw[hyperedgeB] (p1) to[out=+45, in=180-80, relative,looseness=3.5] (p2)
                         to[out=-10, in=180+10, relative] (p3)
                         to[out=+35, in=180-30, relative] (p1);                         
}

\foreach \i in {0,...,4}{
  \pgfmathtruncatemacro\j{\i+1}
  \node at (c\i) {$\frac{2}7$};
  \node at ({90+72*\i}:{\R+\rc+0.35}) {$V_{\j}$};
}
\node at (0,{-\R-\rc-0.9}) {$H_n$};
\end{tikzpicture}
}
\hskip 2em
\tikzset{
vtxBig/.style={inner sep=12pt, outer sep=0pt, circle, fill=white,draw}, 
}
\vc{
\begin{tikzpicture}[scale=1.4]
\def\RC{1.5}                    
\coordinate (a) at (-3.3, 1);   
\coordinate (b) at (-3.3,-1);   
\def\rAB{0.6}                   
\def\d{0.35}                    
\def\dC{1.15}                   
\def\rin{0.3}  

\draw (0,0) circle (\RC);
\draw (a) circle (\rAB);
\draw (b) circle (\rAB);

\begin{scope}[rotate=-30]
\draw (30:0.9) coordinate(x1) node[vtxBig]{};
\draw (150:0.9) coordinate(x2) node[vtxBig]{};
\draw (270:0.9) coordinate(x3) node[vtxBig]{};
\draw
(30:0.8) coordinate(1)
(150:0.8) coordinate(2)
(270:0.8) coordinate(3)
;
\draw[hyperedge] (1) to[out=210,in=330] (2) to[out=330,in=90] (3) to[out=90,in=210] (1);
\foreach \ashift in {30,150,270}{
\draw
(\ashift:0.8)++(\ashift+60:0.1) coordinate(1)
++(\ashift+60:0.3) coordinate(2)
(\ashift+120:0.8)++(\ashift+120+90:-0.2) coordinate(3)
;
\draw[hyperedge] (1) to[bend left] (2) to[bend left=5] (3) to[bend left=5] (1);
}


\end{scope}
\path (a) +(200:\d) coordinate (aAAB1) +(280:\d) coordinate (aAAB2)
          +(320:\d) coordinate (aABC)
          +(-20:\d)  coordinate (aAAC1) +(60:\d)  coordinate (aAAC2);
\path (b) +(130:\d) coordinate (bAAB)
          +(60:\d)  coordinate (bABC)
          +(20:\d)   coordinate (bBBC1) +(300:\d) coordinate (bBBC2);
\path (160:\dC) coordinate (cAAC)
      (180:\dC) coordinate (cABC)
      (200:\dC) coordinate (cBBC);

\draw[hyperedge] (aAAB1) to[bend left=5] (bAAB) to[bend left=5] (aAAB2) to[bend left] (aAAB1);
\draw[hyperedge] (aABC) to[bend left=30] (bABC) to[bend left=15] (cABC) to[bend left=15] (aABC);
\draw[hyperedge] (aAAC2) to[bend left] (aAAC1) to[bend left=5] (cAAC) to[bend left=5] (aAAC2);
\draw[hyperedge] (bBBC2) to[bend left=5] (cBBC) to[bend left=5] (bBBC1) to[bend left] (bBBC2);

\draw (60:1.75) node{$C$};
\path (a) ++(180:0.85) node{$A$};
\path (b) ++(180:0.85) node{$B$};
\draw (-0.9,-2.1) node{$X_n$};
\end{tikzpicture}
}
\end{center}
\caption{$3$-graphs $H_n$ and $X_n$.}\label{fig:constructions}
\end{figure}
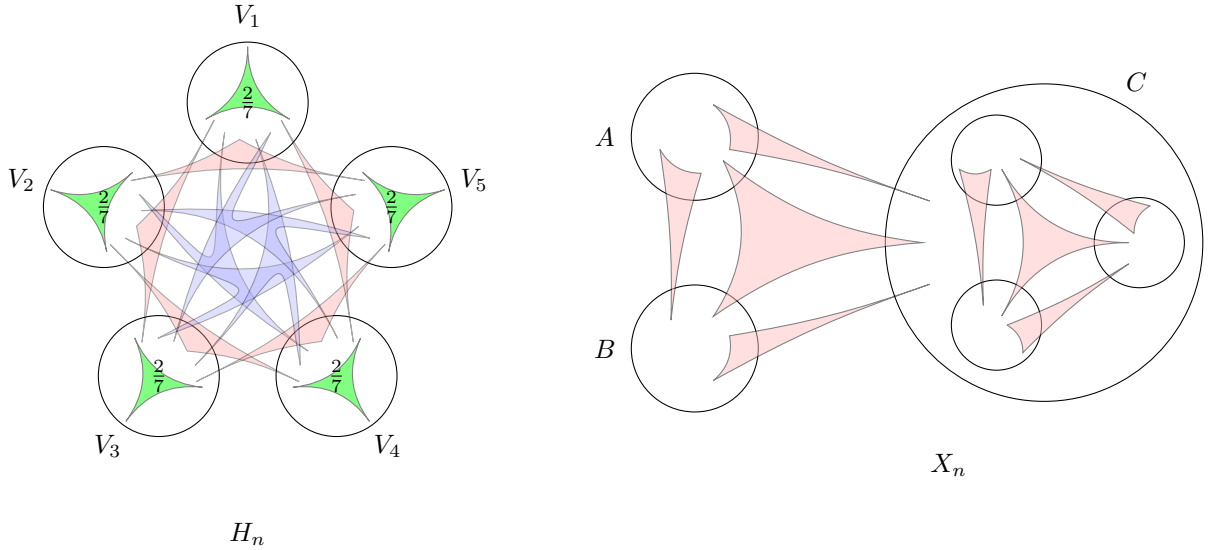

Let $H_n$ be a 3-graph on $N=30n$ vertices, which are partitioned into five parts $V_1,\ldots,V_5$ of equal sizes.
Include every triple meeting three distinct parts, no triple meeting exactly two parts, and add
the $K_4^-$-free 3-graph of density $2/7$ described in the next paragraph inside each part; See the left image in Figure~\ref{fig:constructions} for an illustration.

Recall that $\pi\left(K_4^-\right)$ and $\pi_{\ell_2}\left(K_4^-\right)$ are conjectured to 
be $2/7$ and $4/43$, respectively.
In fact, it is conjectured that both variants have the same unique maximizer, namely
the sequence of iterated blow-ups of the so-called $(6,3,2)$-design
\[
T=\{124,126,134,135,156,235,236,245,346,456\}.
\]
For more details, see~\cite{MR4421399,survey}.
Consequently, the limit of the triple density of $H_n$ is equal to
\[
\frac{4\cdot3}{5^2}+\frac1{5^2}\cdot\frac27=\frac{86}{175},
\]
and analogously the limit of the normalized $\textrm{co}_2\left(H_n\right)$ is
\[
\frac45\left(\frac35\right)^2+\frac1{5^3}\cdot\frac4{43}
=\frac{1552}{5375}.
\]

It remains to prove that $H_n$ is $(F_{4^-,2})$-free.
Suppose, for the sake of contradiction, that there are six vertices $1,\ldots,6$ in $H_n$ forming a copy of $(F_{4^-,2})$.
If $5$ and $6$ were from the same part of $H_n$, then the four triples $x56$ force all of $1$, $2$, $3$ and $4$ into one part (in fact, the same one), contradicting the $K_4^-$-freeness of the parts.
Thus $5$ and $6$ must be in distinct parts, say $V_1$ and $V_2$.
This time, each triple $x56$ forces $x \in \{1,2,3,4\}$ to avoid $V_1 \cup V_2$.
By the pigeonhole principle, at least two of the four vertices must be in the same part.
However, the presence of the triples $123$, $124$ and $134$ again yields
that all the vertices $1$, $2$, $3$ and $4$ must be in the same part, which contradicts the $K_4^-$-freeness of the parts one more time.

\subsection{A construction for $\pi\left(F_{4,2}\right)\ge 0.618782$ and $\ell_2\left(F_{4,2}\right) \ge 0.4240915$}
\label{sec:constructions:F42}

This time, we construct $n$-vertex 3-graphs $X_n$ with their vertex-sets partitioned into three parts $A_n$, $B_n$ and $C_n$. 
$X_n$ contains all the triples that has exactly one vertex in $C_n$, all the triples with two vertices in $A_n$ and one in $B_n$, and finally $C_n$ induces the $3$-graph that is asymptotically extremal for $\pi_{\ell_2}(K_4)$, which is known to have density $5/9$.
An illustration is depicted on the right in Figure~\ref{fig:constructions}.

It is routine to check that $X_n$ are $F_{4,2}$-free.
Moreover, optimizing the sizes of the parts in the above construction yields that $\pi\left(F_{4,2}\right)\geq \frac12+\frac{2\sqrt{14}}{63} \sim 0.618782$. Note that the respective relative sizes of the $A_n$, $B_n$ and $C_n$ are equal to  $\left(\frac13+\frac{\sqrt{14}}{21},\frac16+\frac{\sqrt{14}}{42},\frac12-\frac{\sqrt{14}}{14}\right) \sim (0.511507495,\ 0.255753748,\ 0.232738758)$.

On the other hand, setting the relative sizes of $A_n$, $B_n$ and $C_n$ to be $(0.543709,\ 0.246074,\ 0.210217)$, which corresponds to a numerical solution of a degree-6 polynomial related to the $\textrm{co}_2$-optimization of $X_n$,
yields that $\pi_{\ell_2}\left(F_{4,2}\right) \geq 0.4240915$.

We note that Claude claimed it is possible to further improve the construction for Tur\'an density by splitting $A$ into more classes. 
However, we opted to present an inferior but nicer construction,
as we do not expect the aforementioned improvement to yield the best possible construction.

\subsection{A construction for $\pi\left(J_5\right)\ge \frac9{14}$ and $\ell_2\left(J_5\right) \ge \frac{27}{61}$}
\label{sec:constructions:J5}

Recall that $PG(2,3)$ is a projective plane with 13 points, 13 lines, every line contains 4 points, every point lies on 4 lines, and every pair of points lies on a unique line.
Let $P$ be the $13$-vertex 3-graph defined on the points of the $PG(2,3)$, where the triples of $P$ are exactly the triples of points of $PG(2,3)$ that are not collinear.
Note that $P$ has 
\[
\binom{13}{3}-13\binom43
=286-52
=234
\]
triples.

We claim $P$ is $J_5$-free. Indeed, the link of every vertex of $P$ is isomorphic to $K_{3,3,3,3}$, which is $K_5$-free. However, the link graph of the first vertex in $J_5$ is $K_5$.
A moment of thought reveals that the link graph of any vertex of the iterated blow-up of $P$
is a disjoint union of complete $4$-partite graphs.
In particular, its vertices also have $K_5$-free links, which in turn yields that $\pi(J_5) \geq 9/14$ and $\pi_{\ell_2}(J_5) \geq 27/61$.

\section{New results using CSDP}

While tackling different extremal problems using the method described in Section~\ref{sec:method}, we have noticed that, compared to a decade ago, CSDP can now handle higher levels of the flag algebra hierarchy for several important open problems.
The one is the Tur\'an density of a 3-graph $C_5$. The reported CSDP allocation is $1{,}855$ GB, allowing this instance to run on a machine with approximately $2$ TB of RAM, subject to additional memory overhead.
We obtained using CSDP that $\pi(C_5) \leq  0.465031$.
Note that our method gave a slightly worse upper bound $\pi(C_5) \leq  0.46507$.

In the rest of this section, we discuss our new results for undirected graphs.
Let us recall that the \emph{inducibility} of a fixed graph $H$ is defined as
\[
I(H) := \lim_{n\to\infty} \max_{G \in \mathcal{F}_n} P(H,G). 
\]
This notion was introduced by Pippenger and Golumbic~\cite{MR401552}, who also observed that the limit always exist.
The inducibility of every graph on 4 vertices is known~\cite{EvenZohar2014} with the sole exception of $P_4$, the path on 4 vertices.
Even-Zohar and Linial~\cite{EvenZohar2014}
gave a construction showing 
$I(P_4) \geq 1173/5824 \approx 0.2014$.
The best known upper bound on $I(P_4)$ was by Vaughan~\cite{flagmatic} who used flag algebras on 8 vertices.
Recall that in the setting of graphs, it holds that $f_8=12{,}346$ and $f_9 = 274{,}668$.
However, since $P_4$ is self-complementary, we can exploit the symmetry and the space to $137{,}352$ in a well-known fashion (see, e.g., Kiem, Pokutta and Spiegel~\cite{kiem2023fourcolorramseymultiplicitytriangles}).
Calculating with $9$-vertex subgraph densities yields the following improved bound. 

\begin{theorem}
    $I(P_4) \leq 0.20332375$.
\end{theorem}

The solutions for $\pi(C_5)$ and $I(P_4)$ are available together with the other calculations.

\section{Graph drawing and other applications}
\def\ucr#1{{\text{\rm cr}(#1)}}
\def\rcr#1{\overline{{\text{\rm cr}}}(#1)}
Given a graph $G$, recall that the {\em crossing number} $\ucr{G}$ of $G$ is the minimum number of pairwise crossings of edges in a drawing of $G$ in the plane.
Zarankiewicz conjectured that the crossing number  $\ucr{K_{m,n}}$ of the complete bipartite graph with the parts of respective sizes $m$ and $n$ is equal to
\[
Z(m,n):=\bigg\lfloor{\frac{m}{2}}\bigg\rfloor \bigg\lfloor{\frac{m-1}{2}}\bigg\rfloor \bigg\lfloor{\frac{n}{2}}\bigg\rfloor \bigg\lfloor{\frac{n-1}{2}}\bigg\rfloor.
\]
A closely related conjecture is the one of Hill that states that the crossing number $\ucr{K_{n}}$ of the complete graph $K_n$ is
\[
H(n):=
\frac14 
\left\lfloor \frac{n}{2} \right\rfloor
\left\lfloor \frac{n-1}{2} \right\rfloor
\left\lfloor \frac{n-2}{2} \right\rfloor
\left\lfloor \frac{n-3}{2} \right\rfloor
\]
for all $n \geq 3$.
Both conjectures were shown to be true by OpenAI~\cite{OAI:The-crossing-number-of-complete-bipartite-graphs-September-23-2026,OAI:The-crossing-number-of-complete-graphs-September-23-2026}.
Before these proofs, the best lower bounds were
$\ucr{K_{n,n}} \geq \left(0.9118-o(1)\right)\cdot Z(n,n)$
and
 $\ucr{K_{n}} \ge \left(0.985-o(1)\right)\cdot H(n)$ using flag algebras~\cite{BaloghHill2019,Balogh2023}.
 By utilizing the solver and the method from this paper,
Balogh, Layman, Lidick\'{y}, Pfender, Salazar, Spiro and Volec~\cite{sadpaper} showed
$\ucr{K_{n,n}} \geq \left(0.96937-o(1)\right) \cdot Z(n,n)$
and
$\ucr{K_{n}} \geq \left(0.9899-o(1)\right)\cdot H(n)$.

Balogh, Clemen and Lidick\'y~\cite{MR4421399} conjectured that $\pi_{\ell_2}(K_5^{(4)}) = \frac{31}{64}$, where $K_5^{(4)}$ is the complete 4-graph with $5$ vertices.
In a follow-up paper with Pfender and Sankar~\cite{4uniform} we tackle Turán densities of various $4$-graphs, and, in particular, prove this conjecture. 
We were notified by Corrêa~\cite{zenodo} that he also resolved the conjecture.

Erd\H{o}s~\cite{Erdos84} was interested in the Tur\'an density of $C_4$ in hypercubes denoted by $\pi_{\mathcal{Q}}(C_4)$.
He offered \$100 for for the question if $\pi_{\mathcal{Q}}(C_4)=1/2$.
The upper bound $\pi_{\mathcal{Q}}(C_4) \leq 0.62284$ was obtained by Chung~\cite{Chung1992}.
It was improved by Thomason and Wagner~\cite{Thomason2009} by a computer assisted proof to 0.62256.
Balogh, Hu, Lidick\'y and Liu~\cite{Balogh2014Cube} used flag algebras to obtain $\pi_{\mathcal{Q}}(C_4) \leq 0.6068$.
Solving a flag algebra relaxation via the method from this paper Layman and Lidick\'y~\cite{hypercube} claimed the upper bound can be pushed further below 0.59143.

\section{Conclusion}
We demonstrated a systematic approach to solve flag algebras problems using semidefinite relaxations that works for larger instances than the previous methods. 
We note that the bounds presented in this paper could be improved by giving the solver more time as well as by including so-called \emph{regularity constraints}, which correspond to a flag algebra translation of Zykov's symmetrization argument~\cite{Zykov}.

We expect that solving general relaxations with the 7-vertex subgraph densities is as far as flag algebras can be pushed on 3-graphs with the current methods, given the growth of the number of non-isomorphic configurations (recall Table~\ref{tab:counts}). Indeed, $f_8 > 10^{12}$ makes approaching 8-vertex densities in a nearby future simply by waiting for stronger computers very unlikely.
In the setting of undirected graphs, however, it looks promising that calculations on 10-vertex densities might be within reach.

\section{Use of AI}
The idea of tackling hypergraph Turán problems with this method, our implementation of the solver, the bound verificators, as well as the particular applications were all done solely by the authors. 
ChatGPT-6 was used to find the constructions presented in Sections~\ref{sec:constructions:F4-2}-\ref{sec:constructions:J5}. 
The entries with a dagger $\dagger$ in Table~\ref{tab:results} had exact rounding obtained using ChatGPT-6 and the result was verified by a program written by the authors of this paper.
It was also used to correct English and look for mathematical errors in the manuscript.
The authors have checked the correctness of all the content in this paper.

\section{Acknowledgment}
This material is based upon work supported by the National Science Foundation under Grant No. DMS-1928930, while the authors were in residence at the Simons Laufer Mathematical Sciences Institute in Berkeley, California, during the spring of 2025.
The calculations were performed mainly on the Alderaan cluster.
The Alderaan cluster was funded by the National Science Foundation Campus Computing program grant OAC-2019089 and a contribution from the Simons Foundation.
It is operated by the Center for Computational Mathematics, Department of Mathematical and Statistical Sciences, College of Liberal Arts and Sciences at University of Colorado, Denver.
The work of both authors was also partially supported by the MEYS of the Czech Republic under the INTER-EXCELLENCE II program (project No. LUAUS26294).

We thank Corrêa for pointing out a mistake in calculation of lower bound $\pi_{\ell_2}(J_4)$ in the survey paper~\cite{MR4421399}.

\bibliographystyle{plainurl}
\bibliography{rev}

\end{document}